\documentclass[final,3p,10pt,authoryear]{elsarticle}

\makeatletter
\def\ps@pprintTitle{%
	\let\@oddhead\@empty
	\let\@evenhead\@empty
	\let\@oddfoot\@empty
	\let\@evenfoot\@oddfoot
}
\makeatother

\usepackage{graphicx} 

\usepackage{bm}
\usepackage{amsmath}
\usepackage{amssymb}
\usepackage{IEEEtrantools}

\newcommand{\mat}[1]{\mathbf{#1}}
\newcommand{\argmin}{\operatorname*{arg\,min}}

\newcommand{\tr}[1]{#1^{\mathsf{T}}}
\newcommand{\inv}[1]{#1^{-1}}

\newcommand{\diag}[1]{\mathbf{diag}{#1}}

\usepackage{color}

\usepackage{algorithm,algpseudocode}
\algblockdefx[NAME]{Algorithm}{EndAlgorithm}%
[2][Unknown]{{\bf Algorithm}: #1(#2)}%
{{\bf End Algorithm}}
\algnewcommand\algorithmicto{\textbf{to}}
\algrenewtext{For}[3]%
{\algorithmicfor\ #1 $\gets$ #2 \algorithmicto\ #3 \algorithmicdo}
\algnewcommand\algorithmicReturn{\textbf{return}\;\;}
\algnewcommand\algorithmicBreak{\textbf{break}\;\;}
\algnewcommand\algorithmicinput{\textbf{Input:}}
\algnewcommand\Input{\item[\algorithmicinput]}
\algnewcommand\algorithmicoutput{\textbf{Output:}}
\algnewcommand\Output{\item[\algorithmicoutput]}
\algnewcommand\algorithmicName[1]{\textbf{Algorithm:} #1}
\algnewcommand\Name{\item[\algorithmicName]}

\makeatletter
\long\def\pprintMaketitle{\clearpage
	\iflongmktitle\if@twocolumn\let\columnwidth=\textwidth\fi\fi
	\resetTitleCounters
	\def\baselinestretch{1}%
	\printFirstPageNotes
	\begin{center}%
		\thispagestyle{pprintTitle}%
		\def\baselinestretch{1}%
		\Large\@title\par\vskip18pt
		\normalsize\elsauthors\par\vskip10pt
		\footnotesize\itshape\elsaddress\par\vskip36pt
		\ifvoid\absbox\else\unvbox\absbox\par\vskip10pt\fi
		\ifvoid\keybox\else\unvbox\keybox\par\vskip10pt\fi
	\end{center}%
	\gdef\thefootnote{\arabic{footnote}}%
}
\makeatother

\begin{document}
	\begin{frontmatter}
		\title{A note on alternating direction method of multipliers with generalized augmented terms for constrained sparse least absolute deviation}
		\author[umass]{Yuki Itoh}
		
		\author[umass]{Mario Parente}
		\address[umass]{Department of Electrical and Computer Engineering, University of Massachusetts, Amherst}
		\begin{abstract}
			This technical note is an ancillary material for our research paper~\citep{Itoh2019}. We discuss an alternating direction method of multipliers with generalized augmented terms (ADMM-GAT) and introduce a generalized residual balance technique for efficiently employing ADMM-GAT. These techniques are applied to least absolute deviation and its constrained version and their algorithmic details are presented. These algorithms are used for the implementation of the method described in~\citep{Itoh2019}.
		\end{abstract}
	\end{frontmatter}

The alternating direction method of multipliers (ADMM)~\citep{Boyd2010} is widely used in solving convex optimization problems. 
ADMM with generalized augmented terms (ADMM-GAT) is briefly mentioned in Section 3.4.2 in the tutorial paper~\citep{Boyd2010} on ADMM, which says that it can be cast as a standard ADMM by imposing an additional equality constraint. However, with this formulation, spectral penalty parameters in the augmented terms are considered to be constant, which hinders the automatic adjustment of the magnitude of the penalty terms, such as residual balancing, during the optimization. In order to take full advantage of the generalization, a technique that performs automatic adjustment of the penalty parameters for ADMM-GAT is necessary. We provide a new formulation of residual balancing for ADMM-GAT to further speed up the convergence of ADMM-GAT. 

Section~\ref{sec:admm-gat} describes the theory of ADMM-GAT and the new residual balancing technique. Section~\ref{sec:example_admm-gat} presents the application of ADMM-GAT to least absolute deviation (LAD) and constraint basis pursuit (CBP) and Their pseudo codes are given. We also provide a compromised version of the residual balancing for batch processing of these two applications. Finally, Section~\ref{sec:cslad} shows how to convert constrained sparse LAD (CSLAD) to a CBP problem.

\section{ADMM-GAT}
\label{sec:admm-gat}
We consider a general convex optimization problem for which ADMM can be used:
\begin{IEEEeqnarray*}[]{l"c}
	\begin{IEEEeqnarraybox}[][t]{l'l}
		\underset{\bm{x},\bm{y}}{\text{minimize}} & f(\bm{x})+g(\bm{z}) \\
		\text{subject to} & \mat{A}\bm{x} + \mat{B}\bm{z} = \bm{c}.
	\end{IEEEeqnarraybox}\IEEEyesnumber\label{cvx:general}
\end{IEEEeqnarray*}
where $\bm{x} \in \mathbb{R}^{m}$, $\bm{z} \in \mathbb{R}^n$  $\bm{c} \in \mathbb{R}^p$, $\mat{A} \in \mathbb{R}^{p\times m}$, $\mat{B} \in \mathbb{R}^{p\times n}$, and $f$ and $g$ are convex functions with respect to $\bm{x}$ and $\bm{y}$, respectively. Let us consider a Lagrangian with a general augmentation term:
\begin{IEEEeqnarray*}[]{r,c,l}
	\mathcal{L}(\bm{x},\bm{z},\bm{d}) &=& f(\bm{x})+g(\bm{z}) + \tr{\bm{y}}(\mat{A}\bm{x} + \mat{B}\bm{z} - \bm{c}) + \frac{1}{2}{\|\mat{F}(\mat{A}\bm{x} + \mat{B}\bm{z} - \bm{c})\|}_2^2
\end{IEEEeqnarray*}
where $\bm{y} \in \mathbb{R}^p$ is a vector of Lagrangian multipliers and $\mat{F} \in \mathbb{R}^{q\times p}$ is a general spectral penalty parameter.  Let $\mat{P}= \tr{{\mat{F}}}\mat{F}$, then the Lagrangian is expressed only with $\mat{P}$:
\begin{IEEEeqnarray*}[]{r,c,l}
	\mathcal{L}(\bm{x},\bm{z},\bm{d}) &=& f(\bm{x})+g(\bm{z}) + \tr{\bm{y}}(\mat{A}\bm{x} + \mat{B}\bm{z} - \bm{c}) + \frac{1}{2}\tr{(\mat{A}\bm{x} + \mat{B}\bm{z} - \bm{c})}\mat{P}(\mat{A}\bm{x} + \mat{B}\bm{z} - \bm{c})
\end{IEEEeqnarray*}

The scheme of the ADMM-GAT is same as that of ADMM, which is the iteration of the following problems:
\begin{IEEEeqnarray*}[]{r,c,l}
	\bm{x}^{(k+1)} &\gets& \argmin_{\bm{x}} \mathcal{L}(\bm{x},\bm{z}^{(k)},\bm{y}^{(k)}) \IEEEyesnumber\IEEEyessubnumber\label{eq:gadmm_update_x}\\
	\bm{z}^{(k+1)} &\gets& \argmin_{\bm{z}} \mathcal{L}(\bm{x}^{(k+1)},\bm{z},\bm{y}^{(k)}) \IEEEyessubnumber\label{eq:gadmm_update_z}\\
	\bm{y}^{(k+1)} &\gets& \bm{y}^{(k)} + \mat{P} \bigl({\mat{A}\bm{x}^{(k+1)} -\mat{B}\bm{z}^{(k+1)}-\bm{c}}\bigr), \IEEEyessubnumber\label{eq:gadmm_update_d}
\end{IEEEeqnarray*}
where the superscript $(k)$ indicates the number of iteration. The equation~\eqref{eq:gadmm_update_d} is a dual ascent step. Similarly, the scaled form of the augmented Lagrangian is 
\begin{IEEEeqnarray*}[]{r,c,l}
	\mathcal{L}(\bm{x},\bm{z},\bm{d}) &=& f(\bm{x})+g(\bm{z}) + \tr{\bm{d}}\tr{{\mat{F}}}\mat{F}(\mat{A}\bm{x} + \mat{B}\bm{z} - \bm{c}) + \frac{1}{2}{\|\mat{F}(\mat{A}\bm{x} + \mat{B}\bm{z} - \bm{c})\|}_2^2 \\
	&=& f(\bm{x})+g(\bm{z}) + \tr{\bm{d}}\mat{P}(\mat{A}\bm{x} + \mat{B}\bm{z} - \bm{c}) + \frac{1}{2}\tr{(\mat{A}\bm{x} + \mat{B}\bm{z} - \bm{c})}\mat{P}(\mat{A}\bm{x} + \mat{B}\bm{z} - \bm{c})
\end{IEEEeqnarray*}
where $\bm{d} = \inv{\mat{P}}\bm{y}$ is the vector of a scaled Lagrangian multipliers and its update scheme is
\begin{IEEEeqnarray*}[]{r,c,l}
	\bm{x}^{(k+1)} &\gets& \argmin_{\bm{x}} \mathcal{L}(\bm{x},\bm{z}^{(k)},\bm{d}^{(k)}) \IEEEyesnumber\IEEEyessubnumber\label{eq:gadmm_update_x_scaled}\\
	\bm{z}^{(k+1)} &\gets& \argmin_{\bm{z}} \mathcal{L}(\bm{x}^{(k+1)},\bm{z},\bm{d}^{(k)}) \IEEEyessubnumber\label{eq:gadmm_update_z_scaled}\\
	\bm{d}^{(k+1)} &\gets& \bm{d}^{(k)} + \bigl({\mat{A}\bm{x}^{(k+1)} -\mat{B}\bm{z}^{(k+1)}-\bm{c}}\bigr), \IEEEyessubnumber\label{eq:gadmm_update_d_scaled}
\end{IEEEeqnarray*}

\subsection{Residual-balancing for ADMM-GAT}
\label{sec:gadmm_update_rhoP}
Residual balancing is a common heuristic used for the automatic adjustment of spectral penalty parameters. It considers the primal residual:
\begin{IEEEeqnarray*}[]{r,c,l}
	\bm{r}^{(k+1)} = \mat{A}\bm{x}^{(k+1)} + \mat{B}\bm{z}^{(k+1)} - \bm{c},
\end{IEEEeqnarray*}
and dual residual:
\begin{IEEEeqnarray*}[]{r,c,l}
	\bm{s}^{(k+1)} &=& \tr{\mat{A}}\mat{P}\mat{B}(\bm{z}^{(k+1)}-\bm{z}^{(k)})
\end{IEEEeqnarray*}
and keeps these values within a same order of magnitude. Unlike a normal ADMM with a scalar spectral penalty parameter, the adjustment of the values of $\mat{P}$ is not straightforward. It is necessary to evaluate how each component of the primal residual associated with each element of $\mat{P}$ differs from the counterpart of the dual residual. This can be evaluated when $\mat{P}$ is a diagonal matrix. Let us define $\diag{(\mat{P})} = [P_1,P_2,\ldots,P_p]$. The component of the primal residual associated with $P_l$ is the $l$th element $r_l^{(k+1)}$ of $\bm{r}^{(k+1)}$. For the dual residual component associated with $P_l$, we consider an expansion:
\begin{IEEEeqnarray*}[]{r,c,l}
	\tr{\bm{1}}\bm{s}^{(k+1)} &=&  \sum_{l=1}^{p} P_l \bigl({\bm{b}}^l\cdot(\bm{z}^{(k+1)}-\bm{z}^{(k)})\bigr) \bigl({\bm{a}}^l\cdot\bm{1}_m\bigr),
\end{IEEEeqnarray*}
where ${\bm{a}}^l$ and ${\bm{b}}^l$ are the $l$th row of $\mat{A}$ and $\mat{B}$, respectively. In light of this, we evaluate the magnitude $\bar{s}^{(k+1)}_l$ of the component of the dual residual associated with $P_l$ as
\begin{IEEEeqnarray*}[]{r,c,l}
	\bar{s}^{(k+1)}_l = P_l\cdot\sqrt{\left({{\bigl|\bm{b}^l\bigr|}^2\cdot {\bigl| \bm{z}^{(k+1)}-\bm{z}^{(k)}\bigr|}^2}\right)  \left({{\bigl|\bm{a}^l\bigr|}^2 \cdot \bm{1}_m}\right)},
\end{IEEEeqnarray*}
where ${|\cdot|}^2$ performs the element-wise squares of the vector inside it.

The residual balancing in this case is performed as follows:
$$
P_l^{(k+1)} = 
\begin{cases}
\tau P_l^{(k)} & \text{if\;\;\;} r^{(k)}_l \geq \mu \bar{s}^{(k)}_l \\
\tau / P_l^{(k)} & \text{if\;\;\;} \bar{s}^{(k)}_l \geq \mu r^{(k)}_l \\
P_l^{(k)} & \text{otherwise}
\end{cases}
$$
where $\tau$ and $\mu$ are normally predefined hyper parameters. Typical values are $\tau=10$ and $\mu=2$. 

\section{Example of ADMM-GAT realizations}
\label{sec:example_admm-gat}
We here show the ADMM-GAT of two examples: LAD and CBP. In these examples the spectral penalty parameter matrix $\mat{P}$ is further replaced with $\rho\mat{P}$. This redundant generalization is beneficial when solving the collection of the same problem with partially independent input parameters. With $\rho=1$, we can easily go back to the original ADMM-GAT formulation.
\subsection{ADMM-GAT for LAD}
\label{sec:gadmm_lad}
This section describes a GADMM algorithm to solve Least absolute deviation (LAD):
\begin{IEEEeqnarray*}[]{l"c}
	\textsf{LAD}:
	\begin{IEEEeqnarraybox}[][t]{l'l}
		\underset{\bm{x}}{\text{minimize}} & {\|\bm{h} - \mat{A}\bm{x} \|}_1 
	\end{IEEEeqnarraybox}
\end{IEEEeqnarray*}
where $\bm{x} \in \mathbb{R}^{n}$, $\mat{A} \in \mathbb{R}^{m \times n}$, and $\bm{h} \in \mathbb{R}^m$. Letting $\bm{z} =  \mat{A}\bm{x} - \bm{y}$, the above problem is reformulated as
\begin{IEEEeqnarray*}[]{l"c}
	\begin{IEEEeqnarraybox}[][t]{l'l}
		\underset{\bm{x},\bm{z}}{\text{minimize}} & {\| \bm{z} \|}_1 \\
		\text{subject to} & \mat{A}\bm{x} -\bm{z} = \bm{h}.
	\end{IEEEeqnarraybox}
\end{IEEEeqnarray*}
The scaled version of the generalized augmented Lagrangian of this problem is 
\begin{IEEEeqnarray*}[]{r,c,l}
	\mathcal{L}(\bm{x},\bm{z},\bm{d}) &=& {\|\bm{z}\|}_1 + \rho\tr{\bm{d}}\tr{\mat{F}}\mat{F}(\mat{A}\bm{x}-\bm{z}-\bm{h}) + \frac{\rho}{2}{\|\mat{F}( \mat{A}\bm{x}-\bm{z}-\bm{h})\|}_2^2 \\
	&=& {\|\bm{z}\|}_1 + \frac{\rho}{2}{\|\mat{F}( \mat{A}\bm{x}-\bm{z}-\bm{h}+\bm{d})\|}_2^2 - \frac{\rho}{2}{\|\mat{F}\bm{d}\|}_2^2
\end{IEEEeqnarray*}
where $\rho$ is a scalar spectral penalty parameter, $\mat{F}$ is a matrix of generalized spectral penalty parameters whose inner product matrix, $\tr{\mat{F}}\mat{F} = \mat{P}$, becomes diagonal and $\bm{d}$ is a vector of Lagrangian multipliers. The GADMM algorithm solves the minimization problem by the alternating optimization of the following 
\begin{IEEEeqnarray*}[]{r,c,l}
	\bm{x}^{(k+1)} &\gets& \argmin_{\bm{x}} \frac{\rho}{2}{\left\| \mat{F}\bigl(\mat{A}\bm{x}-\bm{z}^{(k)}-\bm{h} + \bm{d}^{(k)}\bigr) \right\|}_2^2 \\
	\bm{z}^{(k+1)} &\gets& \argmin_{\bm{r}}{\|\bm{z}\|}_1 + \frac{\rho}{2}{\left\| \mat{F}\bigl(\mat{A}\bm{x}^{(k+1)}-\bm{z}-\bm{h} + \bm{d}^{(k)}\bigr) \right\|}_2^2 \\
	\bm{d}^{(k+1)} &\gets& \bm{d}^{(k)} + (\mat{A}\bm{x}^{(k+1)}-\bm{z}^{(k+1)}), \IEEEyesnumber\label{eq:lad_update_d}
\end{IEEEeqnarray*}
where $k$ indicates the number of iteration. The update of $\bm{x}$ is an unconstrained last square problem.
\begin{IEEEeqnarray*}[]{r,c,l}
	\bm{x}^{(k+1)} &\gets& \inv{\bigl(\tr{\mat{A}}\mat{P}\mat{A}\bigr)}\tr{\mat{A}}\mat{P}\bigl( \bm{h} + \bm{z}^{(k)}-\bm{d}^{(k)} \bigr).\IEEEyesnumber\label{eq:lad_update_x}
\end{IEEEeqnarray*}
The update of $\bm{r}$ is only easily defined if $\mat{P}$ is a diagonal matrix so that the minimization with regard to $\bm{z}$ becomes separable for each element. Otherwise, the minimization cannot be done with just one operation. For a diagonal $\mat{P}$, the update equation is expressed as:
\begin{IEEEeqnarray*}[]{r,c,l}
	\bm{z}^{(k+1)} &\gets& \bm{\mathrm{soft}}{\left( \mat{A}\bm{x}^{(k+1)} -\bm{z} + \bm{d}^{(k)},\, \frac{1}{\rho}\cdot \mathrm{\diag}(\inv{\mat{P}}) \right)}.\IEEEyesnumber\label{eq:lad_update_z}
\end{IEEEeqnarray*}
where $\mathrm{soft}(\cdot)$ is a function for performing element-wise soft-thresholding of the vector of the first input (or matrix):
\begin{IEEEeqnarray*}[]{rcl}
	\bm{\mathrm{soft}}(\bm{x},\bm{\kappa}) = \bm{x}_{\bm{\kappa}},
\end{IEEEeqnarray*}
where 
\begin{IEEEeqnarray*}[]{r,c,l}
	\bm{x}_{\bm{\kappa}}[i] &=& 
	\begin{cases}
		0 & \text{if } |\bm{x}[i]| \le \kappa\\
		\mathrm{sign}(\bm{x}[i])\cdot\bigl(|\bm{x}[i]| - \kappa\bigr) & \text{otherwise}.
	\end{cases}
\end{IEEEeqnarray*}
This algorithm converges much faster than the original ADMM especially when the solution of the unconstrained problem is much differ from its constraint version. The drawback is that the matrix inversion in the equation~\eqref{eq:lad_update_x} needs updating whenever $\mat{P}$ is updated.

\subsection{ADMM-GAT for CBP}
\label{sec:gadmm_cbp}
Next we consider a general framework for the constrained basis pursuit de-nosing problem: 
\begin{IEEEeqnarray*}{l}\label{cvx:CBP}
	\begin{IEEEeqnarraybox}[][t]{l'l}
		\begin{IEEEeqnarraybox}[][c]{l'l}
			\underset{\bm{x}}{\text{minimize}}
			& {\bigl{ \|\bm{c}_1 \odot \bm{x}}\bigr\|}_1 \\	
			\text{subject to}
			& \mat{G}\bm{x} = \bm{h} \text{\; and \;}  \bm{x} \succeq \bm{c}_2,
		\end{IEEEeqnarraybox}
	\end{IEEEeqnarraybox}\IEEEyesnumber
\end{IEEEeqnarray*}
where $\bm{x} \in \mathbb{R}^{n}$, $\mat{G} \in \mathbb{R}^{m \times n}$, $\bm{h} \in \mathbb{R}^m$,  $\bm{c}_1 \in \mathbb{R}^n$,  $\bm{c}_2 \in \mathbb{R}^n$, and $\odot$ represents the element-wise multiplication of the two operands. The problem is equivalent to its variable augmented version:
\begin{IEEEeqnarray*}{l}
	\begin{IEEEeqnarraybox}[][t]{l'l}
		\begin{IEEEeqnarraybox}[][c]{l'l}
			\underset{\bm{x},\bm{z}}{\text{minimize}}
			& {\bigl{\|\bm{c}_1 \odot \bm{z}}\bigr\|}_1 \\	
			\text{subject to}
			& \mat{G}\bm{x} = \bm{h} \text{,\;\;}  \bm{z} \succeq \bm{c}_2,  \text{\; and \;} \bm{x}-\bm{z}=\bm{0},
		\end{IEEEeqnarraybox}
	\end{IEEEeqnarraybox}\IEEEyesnumber
\end{IEEEeqnarray*}
which could be solved via alternating minimization. The scaled form of its generalized augmented Lagrangian is defined as
\begin{IEEEeqnarray*}[]{r,c,l}
	\mathcal{L}(\bm{x},\bm{z},\bm{d}) &=&{\bigl\|{ \bm{c}_1 \odot \bm{z}}\bigr\|}_1 + \mathcal{I}_{\bm{z} \succeq \bm{c}_2}(\bm{z}) + \mathcal{I}_{\mat{G}\bm{x} = \bm{h}}(\bm{x}) + \rho\tr{\bm{d}}\tr{\mat{F}}\mat{F}(\bm{x}-\bm{z}) + \frac{\rho}{2}{\|\mat{F}(\bm{x}-\bm{z})\|}_2^2 \\
	&=&{\bigl\|{\bm{c}_1 \odot \bm{t}}\bigr\|}_1 + \mathcal{I}_{\bm{z} \succeq \bm{c}_2}(\bm{z}) + \mathcal{I}_{\mat{G}\bm{x} = \bm{h}}(\bm{x}) + \frac{\rho}{2}{\|\mat{F}(\bm{x}-\bm{z}+\bm{d})\|}_2^2 - \frac{\rho}{2}{\|\mat{F}\bm{d}\|}_2^2 
\end{IEEEeqnarray*}
where $\mathcal{I}_{\bm{z} \succeq \bm{c}_2}(\bm{z})$ is an indicator function of $\bm{z}$ that outputs zero if $\bm{z} \succeq \bm{c}_2$ and $\infty$ otherwise, $\mathcal{I}_{\mat{G}\bm{x}=\bm{h}}(\bm{x})$ is also an indicator one that outputs zero if $\mat{G}\bm{x}=\bm{h}$ and $\infty$ otherwise, $\rho$ is a scalar spectral penalty parameter, $\mat{F}$ is a matrix of generalized spectral penalty parameters whose inner product matrix, $\tr{\mat{F}}\mat{F} = \mat{P}$, becomes diagonal and $\bm{d}\in\mathbb{R}^{L\times 1}$ is a vector of scaled Lagrangian multipliers. Likewise, the minimization is performed via the repetition of three simplified problems:
\begin{IEEEeqnarray*}[]{r,c,l}
	\bm{x}^{(k+1)} &\gets& \argmin_{\bm{x}} L(\bm{x},\bm{z}^{(k)},\bm{d}^{(k)}) \IEEEyesnumber\IEEEyessubnumber\label{eq:cbp_gadmm_x}\\
	\bm{z}^{(k+1)} &\gets& \argmin_{\bm{z}} L(\bm{x}^{(k+1)},\bm{z},\bm{d}^{(k)}) \IEEEyessubnumber\label{eq:cbp_gadmm_z}\\
	\bm{d}^{(k+1)} &\gets& \bm{d}^{(k)} + \bm{x}^{(k+1)}-\bm{z}^{(k+1)}\IEEEyessubnumber\label{eq:cbp_gadmm_d}
\end{IEEEeqnarray*}
where superscripts $(k)$ and $(k+1)$ represent the number of iteration. The last equation~\eqref{eq:cbp_gadmm_d} is a dual-ascent step. Considering the top two problems are formulated as
\begin{IEEEeqnarray*}[]{r,c,l}
	\bm{x}^{(k+1)} &\gets& \argmin_{\bm{x}}\, \frac{\rho}{2}{\left\|\mat{F}\bigl(\bm{x}-\bm{z}^{(k)}+\bm{d}^{(k)}\bigr)\right\|}_2^2 \text{ subject to } \mat{G}\bm{x}=\bm{h} \\
	\bm{z}^{(k+1)} &\gets& \argmin_{\bm{z}}\, {\|\bm{c}_1 \odot \bm{z} \|}_1 + \mathcal{I}_{\bm{z} \succeq \bm{c}_2}(\bm{z}) + \frac{\rho}{2}{\left\|\mat{F}\bigl(\bm{x}^{(k+1)}-\bm{z}+\bm{d}^{(k)}\bigr)\right\|}_2^2,
\end{IEEEeqnarray*}
the first equation~\eqref{eq:cbp_gadmm_x} is analytically solved by
\begin{IEEEeqnarray*}[]{r,c,l}
	\bm{x}^{(k+1)} &\gets& \left( {\mat{I}-\inv{\mat{P}}\tr{\mat{G}}\inv{\bigl(\mat{G}\inv{\mat{P}}\tr{\mat{G}}\bigr)}\mat{G} }\right)\bigl( \bm{z}^{(k)} - \bm{d}^{(k)} \bigr)  + \inv{\mat{P}}\tr{\mat{G}}\inv{\bigl(\mat{G}\inv{\mat{P}}\tr{\mat{G}}\bigr)}\bm{h},\IEEEyesnumber\label{eq:cbp_update_x}
\end{IEEEeqnarray*}
and the equation~\eqref{eq:cbp_gadmm_z} can be also analytically solved by
\begin{IEEEeqnarray*}[]{r,c,l}
	\bm{z}^{(k+1)} &\gets& \bm{\mathrm{soft}}{\left( \bm{\max}{\left(\bm{x}^{(k+1)} + \bm{d}^{(k)},\bm{c}_2\right)},\, \frac{\bm{c}_1}{\rho}\odot \diag{(\inv{\mat{P}})} \right)},\IEEEyesnumber\label{eq:cbp_update_z}
\end{IEEEeqnarray*}
where $\bm{\max}(\cdot)$ is a function for taking element-wise maximum of two vectors (or matrices).

\subsection{Matrix form of CBP and LAD}
Let us consider solving a collection of the problem in the same form. In case of CBP we may have a set $\{\bm{h}\} = \{ \bm{h}_1,\bm{h}_2,\ldots \bm{h}_N\}$ with the other parameters, $\mat{G}$, $\bm{c}_1$, and $\bm{c}_2$, fixed. In case of LAD we may have a set $\{\bm{h}\}$ with the other parameters $\mat{A}$ fixed. 
In this scenario, the LAD problem can be then expressed with a matrix form: 
\begin{IEEEeqnarray*}[]{l"c}
	\begin{IEEEeqnarraybox}[][t]{l'l}
		\underset{\mat{X}}{\text{minimize}} & {\|\mat{H} - \mat{A}\mat{X} \|}_{1,1},
	\end{IEEEeqnarraybox}
\end{IEEEeqnarray*}
where $\mat{H} = \left[{\begin{IEEEeqnarraybox}[][c]{c,c,c,c}\bm{h}_1 & \bm{h}_2 & \ldots & \bm{h}_N\end{IEEEeqnarraybox}}\right] \in \mathbb{R}^{m\times N}$ and ${\|\cdot\|}_{1,1}$ takes the sum of absolute values of all the elements of a matrix. Similarly, CBP is also expressed with a matrix form:
\begin{IEEEeqnarray*}{l}\label{cvx:CBP_mat}
	\begin{IEEEeqnarraybox}[][t]{l'l}
		\begin{IEEEeqnarraybox}[][c]{l'l}
			\underset{\mat{X}}{\text{minimize}}
			& {\bigl{ \|\mat{C}_1 \odot \mat{X}}\bigr\|}_{1,1} \\	
			\text{subject to}
			& \mat{G}\mat{X} = \mat{H} \text{\; and \;}  \mat{X} \succeq \mat{C}_2,
		\end{IEEEeqnarraybox}
	\end{IEEEeqnarraybox}\IEEEyesnumber
\end{IEEEeqnarray*}
where $\mat{C}_1 = [\underbrace{\begin{IEEEeqnarraybox*}[][t]{,c/c/c/c,} \bm{c}_1 & \bm{c}_1 & \dots & \bm{c}_1 \end{IEEEeqnarraybox*}}_{N}]$ and $\mat{C}_2 = [\underbrace{\begin{IEEEeqnarraybox*}[][t]{,c/c/c/c,} \bm{c}_2 & \bm{c}_2 & \dots & \bm{c}_2 \end{IEEEeqnarraybox*}}_{N}]$.

It is possible to separate this problem into each column of $\mat{H}$ and $\mat{X}$, but it would be useful if we could solve this as one problem to avoid redundantly performing matrix inversion whenever the spectral penalty parameters are updated. The redundant formulation of $\mat{P}$ with $\rho\mat{P}$ a compromised solution for this. $\mat{P}$ takes the variation over different row dimensions and $\rho$ does over different columns. We have seen in the previous sections~\ref{sec:gadmm_lad} and \ref{sec:gadmm_cbp} that with the redundant formulation, $\rho$ is taken outside of the matrix inversion. By defining $\rho$ for each column, we could efficiently perform the ADMM-GAT. Let $\rho_i (i=1,2,\ldots,N)$ as $\rho$ for the $i$th column and $\inv{\bm{\rho}} = \left[{\begin{IEEEeqnarraybox}[][c]{c,c,c,c} \inv{\rho_1} & \inv{\rho_2} & \ldots & \inv{\rho_N}\end{IEEEeqnarraybox}}\right] \in \mathbb{R}^{1\times N}$. Then the update equations are straightforwardly obtained. For the CBP problem, the update equations~\eqref{eq:cbp_update_x}, \eqref{eq:cbp_update_z}, and \eqref{eq:cbp_gadmm_d} become
\begin{IEEEeqnarray*}[]{r,c,l}
	\mat{X}^{(k+1)} &\gets& \left( {\mat{I}-\inv{\mat{P}}\tr{\mat{G}}\inv{\bigl(\mat{G}\inv{\mat{P}}\tr{\mat{G}}\bigr)}\mat{G} }\right)\bigl( \mat{Z}^{(k)} - \mat{D}^{(k)} \bigr)  + \inv{\mat{P}}\tr{\mat{G}}\inv{\bigl(\mat{G}\inv{\mat{P}}\tr{\mat{G}}\bigr)}\mat{H},\IEEEyesnumber\IEEEyessubnumber\label{eq:cbp_update_x_mat}\\
	\mat{Z}^{(k+1)} &\gets& \bm{\mathrm{soft}}{\left( \bm{\max}{\left(\mat{X}^{(k+1)} + \mat{D}^{(k)},\mat{C}_2\right)},\, \mat{C}_1 \odot \Bigl({\diag{(\inv{\mat{P}})}\cdot \inv{\bm{\rho}}}\Bigr) \right)},\IEEEyessubnumber\label{eq:cbp_update_z_mat}\\
	\mat{D}^{(k+1)} &\gets& \mat{D}^{(k)} + \mat{X}^{(k+1)}-\mat{Z}^{(k+1)}, \IEEEyessubnumber\label{eq:cbp_update_d_mat}
\end{IEEEeqnarray*}
where $\mat{D}$ is a matrix form of scaled Lagrangian multipliers. 
The update equations are straightforwardly obtained. For the LAD problem, the update equations~\eqref{eq:lad_update_x}, \eqref{eq:lad_update_z}, and \eqref{eq:lad_update_d} becomes
\begin{IEEEeqnarray*}[]{r,c,l}
	\mat{X}^{(k+1)} &\gets& \inv{\bigl(\tr{\mat{A}}\mat{P}\mat{A}\bigr)}\tr{\mat{A}}\mat{P}\bigl( \mat{H} + \mat{Z}^{(k)}-\mat{D}^{(k)} \bigr),\IEEEyesnumber\IEEEyessubnumber\label{eq:lad_update_x_mat}\\
	\mat{Z}^{(k+1)} &\gets& \bm{\mathrm{soft}}{\left( \mat{A}\mat{X}^{(k+1)} -\mat{H} + \mat{D}^{(k)},\, \Bigl({\diag{(\inv{\mat{P}})}\cdot \inv{\bm{\rho}}}\Bigr) \right)},\IEEEyessubnumber\label{eq:lad_update_z_mat}\\
	\mat{D}^{(k+1)} &\gets& \mat{D}^{(k)} + (\mat{A}\mat{X}^{(k+1)}-\mat{Z}^{(k+1)})\IEEEyessubnumber\label{eq:lad_update_d_mat}.
\end{IEEEeqnarray*}

\subsection{Residual balancing for $\rho\mat{P}$ in a matrix form}
In case of matrix form with the redundant formulation of the spectral penalty parameters, the computation of the primal and dual residuals are slightly changed. Here we keep the notation to the general formulation~\eqref{cvx:general}.
The primal residual matrix is
\begin{IEEEeqnarray*}[]{r,c,l}
	\mat{R}^{(k+1)} = \mat{A}\mat{X}^{(k+1)} + \mat{B}\mat{Z}^{(k+1)} - \mat{C},
\end{IEEEeqnarray*}
and the dual residual matrix is 
\begin{IEEEeqnarray*}[]{r,c,l}
	\mat{S}^{(k+1)} &=& \tr{\mat{A}}\mat{P}\mat{B}(\mat{Z}^{(k+1)}-\mat{Z}^{(k)}).
\end{IEEEeqnarray*}
We consider an expansion of the dual residual matrix:
\begin{IEEEeqnarray*}[]{r,c,l}
	\tr{\bm{1}}\mat{S}^{(k+1)}\bm{1}_{N} &=&  \sum_{i=1}^{N} \sum_{l=1}^{p} \rho_iP_l \bigl({\bm{b}}^l\cdot(\bm{z}_i^{(k+1)}-\bm{z}_i^{(k)})\bigr) \bigl({\bm{a}}^l\cdot\bm{1}_m\bigr),
\end{IEEEeqnarray*}
The magnitude of the primal residual associated with $\rho_i$ is
$ r^{(k+1)}_{\rho_i} = {\left\|{\bm{r}_i^{(k+1)}}\right\|}_2$ where $\bm{r}_i^{(k+1)}$ is the $i$th column of $\mat{R}^{(k+1)}$ and that of the dual residual is
\begin{IEEEeqnarray*}[]{r,c,l}
	\bar{s}^{(k+1)}_{\rho_i} = \rho_i\cdot\sqrt{\sum_{l=1}^{p} {|P_l|}^2 \left({{\bigl|\bm{b}^l\bigr|}^2\cdot {\bigl| \bm{z}_i^{(k+1)}-\bm{z}_i^{(k)}\bigr|}^2}\right)  \left({{\bigl|\bm{a}^l\bigr|}^2 \cdot \bm{1}_m}\right)}.
\end{IEEEeqnarray*}
The magnitude of the primal residual associated with $P_l$ is
$ r^{(i+1)}_{P_l} = \sqrt{ \sum_{i=1}^{N}{{ \bigl( {r_{li}^{(k+1)}}  \bigr) }^2}}$ and that of its dual residual is
\begin{IEEEeqnarray*}[]{r,c,l}
	\bar{s}^{(k+1)}_{P_l} = P_l\cdot\sqrt{\left({{\bigl|\bm{b}^l\bigr|}^2\cdot \left({\sum_{i=1}^{N} {\rho_i^2\bigl| \bm{z}_i^{(k+1)}-\bm{z}_i^{(k)}\bigr|}^2}\right)}\right)  \left({{\bigl|\bm{a}^l\bigr|}^2 \cdot \bm{1}_m}\right)},
\end{IEEEeqnarray*}

The residual balancing in this case is performed as follows:
\begin{IEEEeqnarray*}[]{c}
	\rho_i^{(k+1)} = 
	\begin{cases}
		\tau \rho_i^{(k)} & \text{if\;\;\;} r^{(k)}_{\rho_i} \geq \mu \bar{s}^{(k)}_{\rho_i} \\
		\tau / \rho_i^{(k)} & \text{if\;\;\;} \bar{s}^{(k)}_{\rho_i} \geq \mu r^{(k)}_{\rho_i}\\
		\rho_i^{(k)} & \text{otherwise}
	\end{cases}\IEEEyesnumber\label{eq:residualbalancing_gadmm_rho}
\end{IEEEeqnarray*}
and 
\begin{IEEEeqnarray*}[]{c}
	P_l^{(k+1)} = 
	\begin{cases}
		\tau P_l^{(k)} & \text{if\;\;\;} r^{(k)}_{P_l} \geq \mu \bar{s}^{(k)}_{P_l} \\
		\tau / P_l^{(k)} & \text{if\;\;\;} \bar{s}^{(k)}_{P_l} \geq \mu r^{(k)}_{P_l}\\
		P_l^{(k)} & \text{otherwise}
	\end{cases}\IEEEyesnumber\label{eq:residualbalancing_gadmm_P}
\end{IEEEeqnarray*}
Note that we could use different parameters $\tau$ and $\mu$ for the two update rules above. Those two are sequentially performed. It is recommended that the update of the second one uses the updated spectral parameter of the first one. For example, if you perform the update of $\rho_i$ first and $P_l$ second, the updated $\rho_i^{(k+1)}$ will be used for the update of $P_l$. 

\subsection{Pseudo code for CBP and LAD}
\label{sec:pseudo code}
Below are the pseudo codes of CBP and LAD. Practically, the update of the spectral penalty parameters $\rho_i$ and $\mat{P}$ is not performed every iteration. Some of the parameters in the update equations are pre-computed and updated only when $\rho_i$ or $\mat{P}$ is changed. In addition, tolerance is scaled in accordance with the size of the problem.
\begin{algorithm}[!htbp]
	\scriptsize
	\caption{$\mathsf{CBP\,ADMM}\mbox{-}\mathsf{GAT}$($\mat{G},\mat{H}$,$\mat{C}_1$,$\mat{C}_2$,$\epsilon_\mathsf{tol}$,$k_\mathsf{maxiter}$)}\label{alg:cbp_gadmm_pseudocode}
	\begin{algorithmic}[1]
		\Input $\mat{G} \in \mathbb{R}^{m\times n}$, $\mat{H} \in \mathbb{R}^{m\times N}$, $\mat{C}_1 \in \mathbb{R}^{n\times N}$, $\mat{C}_2 \in \mathbb{R}^{n\times N}$
		\Output $\mat{X}^\star  \in \mathbb{R}^{n\times N}$
		\State Set $\rho_i=1 (i=1,\ldots,N)$ and $\mat{P} = \mat{I}$
		\State Pre-compute $\left( {\mat{I}-\inv{\mat{P}}\tr{\mat{G}}\inv{\bigl(\mat{G}\inv{\mat{P}}\tr{\mat{G}}\bigr)}\mat{G} }\right)$, $\inv{\mat{P}}\tr{\mat{G}}\inv{\bigl(\mat{G}\inv{\mat{P}}\tr{\mat{G}}\bigr)}\mat{H}$, and $\mat{C}_1 \odot \Bigl({\diag{(\inv{\mat{P}})}\cdot \inv{\bm{\rho}}}\Bigr) $
		\State Set $R, S\gets\infty$ ($R$ and $S$ are the magnitude of primal and dual residuals, respectively.)
		\State Initializations (if not given):\vspace{-15pt}
		\begin{IEEEeqnarray*}[]{r,c,l}
			\mat{X}^{(0)} &\gets& \inv{\mat{P}}\tr{\mat{G}}\inv{\bigl(\mat{G}\inv{\mat{P}}\tr{\mat{G}}\bigr)}\mat{H}\\
			\mat{Z}^{(0)} &\gets& \bm{\mathrm{soft}}{\left( \bm{\max}{\left(\mat{X}^{(0)},\mat{C}_2\right)},\, \mat{C}_1 \odot \Bigl({\diag{(\inv{\mat{P}})}\cdot \inv{\bm{\rho}}}\Bigr) \right)}\\
			\mat{D}^{(0)} &\gets& \mat{D}^{(k)} + (\mat{X}^{(0)}-\mat{Z}^{(0)}).
		\end{IEEEeqnarray*}\vspace{-20pt}
		\State Set $k=0$ and $\epsilon \gets N\cdot m\cdot \epsilon_\mathsf{tol}$
		\While{$(k<k_\mathsf{maxiter})$ and ($(R>\epsilon)$ or $(S>\epsilon)$)}
		\State Minimize the augmented Lagrangian w.r.t. $\mat{X}$ (Eqn.~\eqref{eq:cbp_update_x_mat}) :
		\begin{IEEEeqnarray*}[]{r,c,l}
			\mat{X}^{(k+1)} &\gets& \left( {\mat{I}-\inv{\mat{P}}\tr{\mat{G}}\inv{\bigl(\mat{G}\inv{\mat{P}}\tr{\mat{G}}\bigr)}\mat{G} }\right)\bigl( \mat{Z}^{(k)} - \mat{D}^{(k)} \bigr)  + \inv{\mat{P}}\tr{\mat{G}}\inv{\bigl(\mat{G}\inv{\mat{P}}\tr{\mat{G}}\bigr)}\mat{H}
		\end{IEEEeqnarray*}
		\State Minimize the augmented Lagrangian w.r.t. $\mat{Z}$ (Eqn.~\eqref{eq:cbp_update_z_mat}) :
		\begin{IEEEeqnarray*}[]{r,c,l}
			\mat{Z}^{(k+1)} &\gets& \bm{\mathrm{soft}}{\left( \bm{\max}{\left(\mat{X}^{(k+1)} + \mat{D}^{(k)},\mat{C}_2\right)},\, \mat{C}_1 \odot \Bigl({\diag{(\inv{\mat{P}})}\cdot \inv{\bm{\rho}}}\Bigr) \right)}
		\end{IEEEeqnarray*}
		\State Dual ascent step (Eqn.~\eqref{eq:cbp_update_d_mat})
		\begin{IEEEeqnarray*}[]{r,c,l}
			\mat{D}^{(k+1)} &\gets& \mat{D}^{(k)} + (\mat{X}^{(k+1)}-\mat{Z}^{(k+1)})
		\end{IEEEeqnarray*}
		\State Update primal residual: $R \gets {\|{\mat{X}^{(k+1)}-\mat{Z}^{(k+1)}}\|}_F$
		\State Update dual residual: $S \gets {\left\|{(\diag{\mat{P}}\cdot\bm{\rho})\odot (\bm{Z}^{(k+1)}-\bm{Z}^{(k)})}\right\|}_F$ 
		\If{$\mod{(k,10)=0}$ or $k=1$}
		\For{$i$}{$1$}{$N$}
		\State $ r^{(k+1)}_{\rho_i} = {\left\|{\bm{r}_i^{(k+1)}}\right\|}_2$ and $\bar{s}^{(k+1)}_{\rho_i} = \rho_i\cdot\sqrt{\sum_{l=1}^{n} {|P_l|}^2 \left({{\bigl| \bm{z}_i^{(k+1)}-\bm{z}_i^{(k)}\bigr|}^2}\right)}$
		\State Update $\rho_i$ by \eqref{eq:residualbalancing_gadmm_rho}
		\EndFor
		\For{$l$}{$1$}{$n$}
		\State$ r^{(k+1)}_{P_l} = \sqrt{ \sum_{i=1}^{N}{{ \bigl( {r_{li}^{(k+1)}}  \bigr) }^2}}$ and $
		\bar{s}^{(k+1)}_{P_l} = P_l\cdot\sqrt{\sum_{i=1}^{N} {\rho_i^2\bigl| \bm{z}_i^{(k+1)}-\bm{z}_i^{(k)}\bigr|}^2}$  (updated $\rho_i$ are used)
		\State Update $P_l$ by \eqref{eq:residualbalancing_gadmm_P}
		\EndFor
		\If{any change in $\bm{\rho}$ or $\mat{P}$}
		\State Update $\left( {\mat{I}-\inv{\mat{P}}\tr{\mat{G}}\inv{\bigl(\mat{G}\inv{\mat{P}}\tr{\mat{G}}\bigr)}\mat{G} }\right)$, $\inv{\mat{P}}\tr{\mat{G}}\inv{\bigl(\mat{G}\inv{\mat{P}}\tr{\mat{G}}\bigr)}\mat{H}$, or $\mat{C}_1 \odot \Bigl({\diag{(\inv{\mat{P}})}\cdot \inv{\bm{\rho}}}\Bigr)$
		\EndIf
		\EndIf
		\State $k\gets k+1$
		\EndWhile
		\State $\mat{X}^{\star} \gets \mat{Z}^{(k)}$
	\end{algorithmic}
\end{algorithm}

\begin{algorithm}[!htbp]
	\scriptsize
	\caption{$\mathsf{LAD\,ADMM}\mbox{-}\mathsf{GAT}$($\mat{A},\mat{H}$,$\epsilon_\mathsf{tol}$,$k_\mathsf{maxiter}$)}\label{alg:lad_gadmm_pseudocode}
	\begin{algorithmic}[1]
		\Input $\mat{A} \in \mathbb{R}^{m\times n}$, $\mat{H} \in \mathbb{R}^{m\times N}$ 
		\Output $\mat{X}^\star \in \mathbb{R}^{n\times N}$
		\State Set $\rho_i=1 (i=1,\ldots,N)$ and $\mat{P} = \mat{I}$
		\State Pre-compute $\inv{\bigl(\tr{\mat{A}}\mat{P}\mat{A}\bigr)}\tr{\mat{A}}\mat{P}$, $\Bigl({\diag{(\inv{\mat{P}})}\cdot \inv{\bm{\rho}}}\Bigr)$
		\State Set $R, S\gets\infty$ ($R$ and $S$ are the magnitude of primal and dual residuals, respectively.)
		\State Initializations (if not given):\vspace{-15pt}
		\begin{IEEEeqnarray*}[]{r,c,l}
			\mat{X}^{(0)} &\gets& \inv{\bigl(\tr{\mat{A}}\mat{P}\mat{A}\bigr)}\tr{\mat{A}}\mat{P}\mat{H}\\
			\mat{Z}^{(0)} &\gets& \bm{\mathrm{soft}}{\left( \mat{A}\mat{X}^{(0)} -\mat{H},\, \Bigl({\diag{(\inv{\mat{P}})}\cdot \inv{\bm{\rho}}}\Bigr) \right)}\\
			\mat{D}^{(0)} &\gets& \mat{D}^{(k)} + (\mat{A}\mat{X}^{(0)}-\mat{Z}^{(0)}).
		\end{IEEEeqnarray*}\vspace{-20pt}
		\State Set $k=0$ and $\epsilon \gets N\cdot n\cdot \epsilon_\mathsf{tol}$
		\While{$(k<k_\mathsf{maxiter})$ and ($(R>\epsilon)$ or $(S>\epsilon)$)}
		\State Minimize the augmented Lagrangian w.r.t. $\mat{X}$ (Eqn.~\eqref{eq:cbp_update_x_mat}) :
		\begin{IEEEeqnarray*}[]{r,c,l}
			\mat{X}^{(k+1)} &\gets& \inv{\bigl(\tr{\mat{A}}\mat{P}\mat{A}\bigr)}\tr{\mat{A}}\mat{P}\bigl( \mat{H} + \mat{Z}^{(k)}-\mat{D}^{(k)} \bigr)
		\end{IEEEeqnarray*}
		\State Minimize the augmented Lagrangian w.r.t. $\mat{Z}$ (Eqn.~\eqref{eq:cbp_update_z_mat}) :
		\begin{IEEEeqnarray*}[]{r,c,l}
			\mat{Z}^{(k+1)} &\gets& \bm{\mathrm{soft}}{\left( \mat{A}\mat{X}^{(k+1)} -\mat{H} + \mat{D}^{(k)},\, \Bigl({\diag{(\inv{\mat{P}})}\cdot \inv{\bm{\rho}}}\Bigr) \right)}
		\end{IEEEeqnarray*}
		\State Dual ascent step (Eqn.~\eqref{eq:cbp_update_d_mat})
		\begin{IEEEeqnarray*}[]{r,c,l}
			\mat{D}^{(k+1)} &\gets& \mat{D}^{(k)} + (\mat{A}\mat{X}^{(k+1)}-\mat{Z}^{(k+1)})
		\end{IEEEeqnarray*}
		\State Update primal residual: $R \gets {\|{\mat{A}\mat{X}^{(k+1)}-\mat{Z}^{(k+1)}}\|}_F$
		\State Update dual residual: $S \gets {\left\|{\tr{\mat{A}}\Bigl( (\diag{\mat{P}}\cdot\bm{\rho}) \odot (\bm{Z}^{(k+1)}-\bm{Z}^{(k)}) \Bigr) }\right\|}_F$ 
		\If{$\mod{(k,10)=0}$ or $k=1$}
		\For{$i$}{$1$}{$N$}
		\State $ r^{(k+1)}_{\rho_i} = {\left\|{\bm{r}_i^{(k+1)}}\right\|}_2$ and $\bar{s}^{(k+1)}_{\rho_i} = \rho_i\cdot\sqrt{\sum_{l=1}^{m} {|P_l|}^2 \left({{\bigl| \bm{z}_i^{(k+1)}-\bm{z}_i^{(k)}\bigr|}^2}\right) \odot\left({{\bigl|\bm{a}^l\bigr|}^2 \cdot \bm{1}_m}\right) }$
		\State Update $\rho_i$ by \eqref{eq:residualbalancing_gadmm_rho}
		\EndFor
		\For{$l$}{$1$}{$m$}
		\State$ r^{(k+1)}_{P_l} = \sqrt{ \sum_{i=1}^{N}{{ \bigl( {r_{li}^{(k+1)}}  \bigr) }^2}}$ and $
		\bar{s}^{(k+1)}_{P_l} = P_l\cdot\sqrt{\sum_{i=1}^{N} {\rho_i^2\bigl| \bm{z}_i^{(k+1)}-\bm{z}_i^{(k)}\bigr|}^2 \odot \left({{\bigl|\bm{a}^l\bigr|}^2 \cdot \bm{1}_m}\right)}$  (updated $\rho_i$ are used)
		\State Update $P_l$ by \eqref{eq:residualbalancing_gadmm_P}
		\EndFor
		\If{any change in $\bm{\rho}$ or $\mat{P}$}
		\State Update $\inv{\bigl(\tr{\mat{A}}\mat{P}\mat{A}\bigr)}\tr{\mat{A}}\mat{P}$ or $\Bigl({\diag{(\inv{\mat{P}})}\cdot \inv{\bm{\rho}}}\Bigr)$
		\EndIf
		\EndIf
		\State $k\gets k+1$
		\EndWhile
		\State $\mat{X}^{\star} \gets \mat{Z}^{(k)}$
	\end{algorithmic}
\end{algorithm}
\newpage
\section{CSLAD}
\label{sec:cslad}
Finally, we consider CSLAD:
\begin{IEEEeqnarray*}[]{l"c}
	\begin{IEEEeqnarraybox}[][t]{l'l}
		\underset{\bm{x}}{\text{minimize}} &  {\| \bm{h}-\mat{G}\bm{x} \|}_1 + {\| \bm{\lambda} \odot \bm{x} \|}_1 \\
		\text{subject to} & \bm{x} \succeq \bm{\gamma}
	\end{IEEEeqnarraybox}
\end{IEEEeqnarray*}
where $\bm{x} \in \mathbb{R}^{n}$, $\bm{h} \in \mathbb{R}^m$, $\mat{G} \in \mathbb{R}^{m\times n}$, $\bm{\lambda} \in \mathbb{R}^{n}$, and $\bm{\gamma} \in \mathbb{R}^n$. We will show that CSLAD comes down to CBP with a variable conversion.
First, letting $\bm{r} = \bm{h} - \mat{G}\bm{x}$, CSLAD is equivalently transformed into:
\begin{IEEEeqnarray*}[]{l"c}
	\begin{IEEEeqnarraybox}[][t]{l'l}
		\underset{\bm{x},\bm{r}}{\text{minimize}} &  {\| \bm{r} \|}_1 + {\| \bm{\lambda} \odot \bm{x} \|}_1 \\
		\text{subject to} & \bm{x} \succeq \bm{\gamma} \text{ and } \bm{r} = \bm{h} - \mat{G}\bm{x}.
	\end{IEEEeqnarraybox}
\end{IEEEeqnarray*}
Then CSLAD is further equivalently converted a general CBP form:
\begin{IEEEeqnarray*}[]{l"c}
	\begin{IEEEeqnarraybox}[][t]{l'l}
		\underset{\bm{u}}{\text{minimize}} &  {\| \hat{\bm{\lambda}} \odot \bm{u} \|}_1  \\
		\text{subject to} & \bm{u} \succeq \hat{\bm{\gamma}} \text{ and } \bm{h} = \hat{\mat{G}}\bm{u}.
	\end{IEEEeqnarraybox}
\end{IEEEeqnarray*}
\begin{IEEEeqnarray*}{t'c}
	where &
	\bm{u} = \left[{\begin{IEEEeqnarraybox}[][c]{c} \bm{x} \\*[-0.5\normalbaselineskip] \bm{r} \end{IEEEeqnarraybox}}\right],\quad
	\hat{\mat{G}} = \left[\begin{IEEEeqnarraybox}[][c]{c,c} \mat{G} & \mat{I}_L\end{IEEEeqnarraybox}\right], \quad
	\hat{\bm{\lambda}} = \left[\begin{IEEEeqnarraybox}[][c]{c} \bm{\lambda} \\*[-0.5\normalbaselineskip] \bm{1}_m \end{IEEEeqnarraybox} \right],\; \text{ and } \quad
	\hat{\bm{\gamma}} =  \left[\begin{IEEEeqnarraybox}[][c]{c} \bm{\gamma}  \\*[-0.3\normalbaselineskip] -\inf \cdot\bm{1}_m\end{IEEEeqnarraybox} \right].
\end{IEEEeqnarray*}
This way the solver of CBP can be used for CSLAD.

In case of a matrix form:
\begin{IEEEeqnarray*}[]{l"c}
	\begin{IEEEeqnarraybox}[][t]{l'l}
		\underset{\mat{X}}{\text{minimize}} &  {\| \mat{H}-\mat{G}\mat{X} \|}_{1,1} + {\| \mat{\Lambda} \odot \mat{X} \|}_{1,1} \\
		\text{subject to} & \mat{X} \succeq \mat{\Gamma}
	\end{IEEEeqnarraybox}
\end{IEEEeqnarray*} 
where $\mat{X} \in \mathbb{R}^{n\times N}$, $\mat{H} \in \mathbb{R}^{m\times N}$, $\mat{\Lambda} \in \mathbb{R}^{n\times N}$, and $\mat{\Gamma} \in \mathbb{R}^{n\times N}$.
Letting $\mat{R} = \mat{H} - \mat{G}\mat{X}$, we have 
\begin{IEEEeqnarray*}[]{l"c}
	\begin{IEEEeqnarraybox}[][t]{l'l}
		\underset{\mat{X},\mat{R}}{\text{minimize}} &  {\| \mat{R} \|}_{1,1} + {\| \mat{\Lambda} \odot \mat{X} \|}_{1,1} \\
		\text{subject to} & \mat{X} \succeq \mat{\Gamma} \text{ and } \mat{R} = \mat{H} - \mat{G}\mat{X}.
	\end{IEEEeqnarraybox}
\end{IEEEeqnarray*}
Then CSLAD is further equivalently converted a matrix form of CBP:
\begin{IEEEeqnarray*}[]{l"c}
	\begin{IEEEeqnarraybox}[][t]{l'l}
		\underset{\mat{U}}{\text{minimize}} &  {\| \hat{\mat{\Lambda}} \odot \mat{U} \|}_{1,1}  \\
		\text{subject to} & \mat{U} \succeq \hat{\mat{\Gamma}} \text{ and } \mat{H} = \hat{\mat{G}}\mat{U}.
	\end{IEEEeqnarraybox}
\end{IEEEeqnarray*}
\begin{IEEEeqnarray*}{t'c}
	where &
	\mat{U} = \left[{\begin{IEEEeqnarraybox}[][c]{,c,} \mat{X} \\*[-0.3\normalbaselineskip] \mat{R} \end{IEEEeqnarraybox}}\right],\quad
	\hat{\mat{G}} = \left[\begin{IEEEeqnarraybox}[][c]{c,c} \mat{G} & \mat{I}_L\end{IEEEeqnarraybox}\right], \quad
	\hat{\mat{\Lambda}} = \left[\begin{IEEEeqnarraybox}[][c]{c} \mat{\Lambda} \\*[-0.5\normalbaselineskip] \bm{1}_{m\times N} \end{IEEEeqnarraybox} \right],\; \text{ and } \quad
	\hat{\mat{\Gamma}} =  \left[\begin{IEEEeqnarraybox}[][c]{c} \mat{\Gamma}  \\*[-0.3\normalbaselineskip] -\inf \cdot\bm{1}_{m\times N}\end{IEEEeqnarraybox} \right].
\end{IEEEeqnarray*}

\bibliographystyle{elsarticle-harv} 

\begin{thebibliography}{2}
\expandafter\ifx\csname natexlab\endcsname\relax\def\natexlab#1{#1}\fi
\providecommand{\url}[1]{\texttt{#1}}
\providecommand{\href}[2]{#2}
\providecommand{\path}[1]{#1}
\providecommand{\DOIprefix}{doi:}
\providecommand{\ArXivprefix}{arXiv:}
\providecommand{\URLprefix}{URL: }
\providecommand{\Pubmedprefix}{pmid:}
\providecommand{\doi}[1]{\href{http://dx.doi.org/#1}{\path{#1}}}
\providecommand{\Pubmed}[1]{\href{pmid:#1}{\path{#1}}}
\providecommand{\bibinfo}[2]{#2}
\ifx\xfnm\relax \def\xfnm[#1]{\unskip,\space#1}\fi
\bibitem[{Boyd et~al.(2010)Boyd, Parikh, Chu, Peleato and Eckstein}]{Boyd2010}
\bibinfo{author}{Boyd, S.}, \bibinfo{author}{Parikh, N.}, \bibinfo{author}{Chu,
  E.}, \bibinfo{author}{Peleato, B.}, \bibinfo{author}{Eckstein, J.},
  \bibinfo{year}{2010}.
\newblock \bibinfo{title}{{Distributed optimization and statistical learning
  via the alternating direcition method of multipliers}}.
\newblock \bibinfo{journal}{Found. Trends Mach. Learn.} \bibinfo{volume}{3},
  \bibinfo{pages}{1--122}.
\bibitem[{Itoh and Parente(2019)}]{Itoh2019}
\bibinfo{author}{Itoh, Y.}, \bibinfo{author}{Parente, M.},
  \bibinfo{year}{2019}.
\newblock \bibinfo{title}{A new method for atmospheric correction and
  de-noising of crism hyperspectral data}.
\newblock \bibinfo{note}{To be submitted}.

\end{thebibliography}

\end{document}